\documentclass[11pt]{amsart}

\usepackage{geometry}
\usepackage{amsfonts, amssymb, amscd}
\numberwithin{equation}{section}

\usepackage[symbol]{footmisc}

\usepackage{todonotes}
 
\usepackage{bm}
\usepackage{verbatim}
\usepackage{amssymb}
\usepackage{mathrsfs}
\usepackage{graphicx}
\usepackage{tikz-cd}
\usepackage{subcaption}
\usepackage{listings}
\usepackage{subfiles}
\usepackage[toc,page]{appendix}
\usepackage{mathtools}
\usepackage{comment}
\usepackage{enumerate}
\usepackage{enumitem}
\usepackage[linesnumbered,ruled]{algorithm2e}
\usepackage[all]{xy}

\usepackage{graphicx}
\graphicspath{{images/}}

\usepackage{appendix}
\usepackage{hyperref}
\newcommand{\bb}{\bm{b}}

\newcommand{\Pp}{\mathbb{P}}
\newcommand{\Qq}{\mathbb{Q}}

\newcommand{\Rr}{\mathbb{R}}

\newcommand{\Zz}{\mathbb{Z}}

\newcommand{\Span}{\operatorname{Span}}

\newcommand{\pet}{\operatorname{pet}}
\newcommand{\PET}{\operatorname{PET}}

\newcommand{\mld}{{\rm{mld}}}

\newcommand{\lct}{\operatorname{lct}}

\newcommand{\Supp}{\operatorname{Supp}}

\newcommand{\mult}{\operatorname{mult}}

\newcommand{\lf}{\lfloor}
\newcommand{\rf}{\rfloor}

\newcommand{\Ii}{{\Gamma}}

\newtheorem{thm}{Theorem}[section]
\newtheorem{conj}[thm]{Conjecture}
\newtheorem{cor}[thm]{Corollary}
\newtheorem{lem}[thm]{Lemma}

\newtheorem{prop}[thm]{Proposition}

\newtheorem{claim}[thm]{Claim}
\theoremstyle{definition}
\newtheorem{rem}[thm]{Remark}

\newtheorem{ex}[thm]{Example}

\theoremstyle{definition}
\newtheorem{defn}[thm]{Definition}

\begin{document}

\title{Uniform rational polytopes for Iitaka dimensions}

\dedicatory{Dedicated to Professor Vyacheslav V. Shokurov on the occasion of his seventieth birthday}

\author{Guodu Chen, Jingjun Han, and Jihao Liu}

\begin{abstract}
In this paper, we continue to develop the theories on functional pairs and uniform rational polytopes. We show that there is a uniform perturbation for Iitaka dimensions of pseudo-effective lc pairs of fixed dimension with DCC coefficients assuming the non-vanishing conjecture. We also show the existence of uniform rational polytopes for Iitaka dimensions of pseudo-effective lc pairs assuming the non-vanishing conjecture. 
\end{abstract}

\address{Institute for Theoretical Sciences, Westlake University, Hangzhou, Zhejiang, 310024, China}
\email{chenguodu@westlake.edu.cn}

\address{Shanghai Center for Mathematical Sciences, Fudan University, Shanghai, 200438, China}
\email{hanjingjun@fudan.edu.cn}

\address{Department of Mathematics, Northwestern University, 2033 Sheridan Rd, Evanston, IL 60208, USA}
\email{jliu@northwestern.edu}

\date{\today}

\maketitle

\tableofcontents

\section{Introduction}
We work over the field of complex numbers $\mathbb C$. 

One of the most important objects in the study of birational geometry is \emph{pairs}, that is, a normal variety $X$ and a boundary $B\geq 0$ on $X$ such that $K_X+B$ is $\Rr$-Cartier. Pairs appear naturally even if one is only interested in studying projective varieties. 

Instead of considering a pair $(X,B:=\sum b_iB_i)$ with fixed coefficients $b_i$, it is natural to study the behavior of singularities when $b_i$'s vary. \cite{HLQ21} introduced and studied (affine) functional pairs, 
$$(X,B(x_1,\ldots,x_c)):=(X,\sum s_i(x_1,\ldots,x_c)B_i),$$ 
where each $s_i(x_1,\ldots,x_c)$ is an affine function of $x_1,\ldots,x_c$, and showed the ACC for log canonical threshold polytopes. In \cite{HLS19}, the last two authors and Shokurov continue to study the theory of functional pairs and developed the theory on \emph{uniform rational polytopes}, a theory which reveals that some birational properties hold for pairs which the coefficients lie in a rational polytope and is independent of the choice of the variety $X$. More precisely, \cite{HLS19} showed the existence of uniform lc (resp. $\Rr$-complementary) rational polytopes \cite[Corollary 5.5]{HLS19}  (resp. \cite[Theorem 5.15]{HLS19}), and finally yields a full solution of Shokurov's boundedness of lc complements conjecture (\cite{Bir19,PS01,PS09,Sho92,Sho00,Sho20}) and the ACC for minimal log discrepancies (mlds) of exceptional singularities. Some embryonic forms of these results on uniform polytopes could be found in \cite{Nak16,Liu18}.
In \cite{CH21}, the first two authors showed the exsitence of uniform rational polytopes for linearity of mlds for surfaces which yields the existence of $(\epsilon, n)$-complements for $(\epsilon, \Rr)$-complementary surface pairs. We refer the readers to \cite{Che20,CX22,HL22,HLL22,HLQ20,Jia20} for other applications of theories on functional pairs and uniform polytopes to ACC conjecture for mlds, the theory of complements, and K-stability theory.

In this paper, we continue to develop the theories on functional pairs and uniform rational polytopes, and apply it to the study on the Iitaka dimensions for pairs. In a sequel paper, we will apply our results to the study on effective Iitaka fibrations.

\medskip

The following theorem is our first main result, see Theorem \ref{thm:dccdecomp1} for the more general statement.

\begin{thm}\label{thm:dccdecomp}
Let $d$ be a positive integer and $\Ii\subset[0,1]$ a DCC set. Assume that Conjecture \ref{conj: non-vanishing} holds for pairs of dimension $d$. Then there exist real numbers $a_1,\dots,a_k\in (0,1]$ and a DCC set $\Ii_0\subset[0,1]$ depending only on $d$ and $\Ii$ such that $\sum_{i=1}^ka_i=1$ and satisfying the following. 

Assume that $(X/Z,B)$ is a $d$-dimensional lc pair such that $B\in\Ii$, all the irreducible components of $B$ are $\Qq$-Cartier, and $K_X+B$ is pseudo-effective over $Z$. Then there exist $\Rr$-divisors $B_1,\dots,B_k\geq 0$ on $X$, such that 
\begin{enumerate}
    \item $\sum_{i=1}^ka_iB_i=B$,
    \item $B_i\in\Ii_0$ for each $i$,
    \item $(X,B_i)$ is lc for each $i$, and
    \item $\kappa(X/Z,K_X+B_i)=\kappa_{\iota}(X/Z,K_X+B)\geq 0$ for each $i$.
\end{enumerate}
\end{thm}

\begin{conj}[Non-vanishing conjecture]\label{conj: non-vanishing}
Let $(X/Z,B)$ be an lc pair. If $K_X+B$ is pseudo-effective over $Z$, then $|(K_X+B)/Z|_{\Rr}\not=\emptyset$. 
\end{conj}

In short, Theorem \ref{thm:dccdecomp} indicates that for pairs $(X,B)$ of fixed dimension with arbitrary DCC coefficients, we can make a uniform decomposition $K_X+B=\sum a_i(K_X+B_i)$, such that the Iitaka dimension of $K_X+B_i$ is equal to the invariant Iitaka dimension of $K_X+B$, and $(X,B_i)$ preserves most nice properties of $(X,B)$. As this decomposition is uniform, we will be able to study the ``effective (invariant) Iitaka fibrations" for these pairs, even if the classical Iitaka fibrations may be undefined for these pairs. We give the following example to illustrate our ideas.

\begin{ex}
Let $\Ii'':=\{\epsilon_i\}_{i=1}^{+\infty}\subset (0,\frac{1}{2})$ be a strictly decreasing sequence of real numbers such that $\lim_{i\rightarrow+\infty}\epsilon_i=0$ (e.g. $\epsilon_i=\frac{1}{i+1}$), and let $\Ii:=\{2-\sqrt{3},\frac{\sqrt{3}}{2},\frac{\sqrt{3}}{2}-\epsilon_i\}_{i\ge1}$ be a DCC set. Let $X:=\mathbb P^1\times\mathbb P^1$, $l_1,l_2,l_3,l_4$ four lines on $X$ of bi-degree $(1,0),(0,1),(1,1),(1,1)$ respectively, such that $l_1,l_2,l_3$ intersect at a unique point $x\in X$ and $x\notin l_4$. Let
$$B_n:=(2-\sqrt{3})l_1+(\frac{\sqrt{3}}{2}-\epsilon_n)l_2+\frac{\sqrt{3}}{2}(l_3+l_4)\in\Ii$$
and $p_1: X\rightarrow\mathbb P^1$ the projection given by $p_1(x_1,x_2)=x_1$. It is clear that for each $n$, $(X,B_n)$ is klt, $\kappa_{\iota}(X,K_X+B_n)=1$, and $p_1$ is the log canonical model of $(X,B_n)$. In other words,  $p_1$ is the \emph{invariant Iitaka fibration} of $K_X+B_n$ (cf. \cite[Definition 2.2]{Hu20}). However, the Iitaka fibration of $K_X+B_n$ is undefined for any $n$ as $\kappa(X,K_X+B_n)=-\infty$.

Nevertheless, we still want to study the effectivity on the invariant Iitaka fibration $p_1$. A natural way is to make a uniform perturbation on the irrational coefficients of $B_n$ to rational coefficients. Unlike the uniform perturbation theory as in \cite{HLS19}, in our scenario, the coefficients of $B_n$ belong to a DCC but not finite set, such perturbation may not exist (indeed, if $\Span_{\mathbb Q}(\{\epsilon_i\}_{i=1}^{+\infty}$) is not finite dimensional over $\Qq$, such perturbation never exists).

On the other hand, we still can make a uniform perturbation on the irrational coefficients of $B_n$ so that $p_1$ is the Iitaka fibration of all induced pairs, even if we cannot guarantee that all induced pairs have rational coefficients. To see this, we consider the functional boundary
$$B_n(t):=(2-t)l_1+(\frac{t}{2}-\epsilon_n)l_2+\frac{t}{2}(l_3+l_4)$$
where $a$ varies. It is clear that $(X,B_n(t))$ is lc for any $t\in (1,2)$. Since $K_X+B_n(t)\sim_{\mathbb R}(\frac{3t}{2}-2-\epsilon_n)l_2$, when $t>\frac{5}{3}>\frac{4+2\epsilon_n}{3}$, $p_1$ is the log canonical model of $K_X+B_n(t)$. In particular, when $t\in (\frac{5}{3},2)\cap\mathbb Q$, $\kappa_{\iota}(X,K_X+B_n(t))=\kappa(X,K_X+B_n(t))=1$ and $p_1$ is the Iitaka fibration of $K_X+B_n(t)$. Therefore, we can discuss the ``effective invariant Iitaka fibration" by considering the effectivity of the Iitaka fibration of $K_X+B_n(t)$ when $t$ varies in $(\frac{5}{3},2)\cap\mathbb Q$.

Interpreting the ideas into the wording as in Theorem \ref{thm:dccdecomp}, we have the following. We let $a_1:=10\sqrt{3}-17, a_2:=18-10\sqrt{3}$, and $\Ii_0:=\{\frac{1}{5},\frac{3}{10},\frac{17}{20}-\epsilon_n,\frac{9}{10}-\epsilon_n\}_{n\ge1}$. Let
$$B_{n,1}:=\frac{1}{5}l_1+(\frac{9}{10}-\epsilon_n)l_2+\frac{9}{10}(l_3+l_4),\text{ and }B_{n,2}:=\frac{3}{10}l_1+(\frac{17}{20}-\epsilon_n)l_2+\frac{17}{20}(l_3+l_4).$$
Then $B_{n,1}=B_n(\frac{9}{5}),B_{n,2}=B_n(\frac{17}{10})$, and we have that $a_1+a_2=1$, $B_{n,1},B_{n,2}\in\Ii_0$, $a_1B_{n,1}+a_2B_{n,2}=B_n,$ $(X,B_{n,i})$ is lc, and $K_X+B_{n,i}$ is pseudo-effective for any $i$. Moreover, $\kappa(X,K_X+B_{n,i})=\kappa_{\iota}(X,K_X+B_n)=1$ and $p_1$ is the Iitaka fibration of $K_X+B_{n,i}$ for any $n$ and $i$.
\end{ex}

The second main theorem of this paper is the existence of a uniform rational polytope for Iitaka dimensions (see Theorem \ref{thm: Uniform kod polytope assume gmm intro1} for the more general statement). Recall that the \emph{rational envelope} of a point $\bm{v}\in\Rr^m$ is the smallest affine subspace containing $\bm{v}$ which is defined over the rationals.

\begin{thm}\label{thm: Uniform kod polytope assume gmm intro}
Let $d$ be a positive integer and $\bm{v}_0:=(v_1^0,\dots,v_m^0)\in\Rr^m$ a point. Assume that Conjecture \ref{conj: non-vanishing} holds for pairs of dimension $d$. Then there exists an open set $V_0\ni\bm{v}_0$ of the rational envelope of $\bm{v}_0$ depending only on $d$ and $\bm{v}_0$ satisfying the following. 

Assume that $X\rightarrow Z$ is a projective morphism between normal quasi-projective varieties, such that $\dim X=d$, and $B_1,\dots,B_m\geq 0$ are Weil divisors on $X$. Let $B(\bm{v}):=\sum_{i=1}^mv_iB_m$ for any point $\bm{v}:=(v_1,\dots,v_m)\in\Rr^m$. Assume that $(X,B(\bm{v}_0))$ is lc, and $K_X+B(\bm{v}_0)$ is pseudo-effective over $Z$. Then for any point $\bm{v}\in V_0$, we have that
\begin{enumerate}
     \item $(X,B(\bm{v}))$ is lc,
    \item $\kappa_{\iota}(X/Z,K_X+B(\bm{v}))=\kappa_{\iota}(X/Z,K_X+B(\bm{v}_0))\geq 0,$ and
    \item if either $\bm{v}\in \Qq^c$ or $\kappa(X/Z,K_X+B(\bm{v}_0))\geq 0$, then 
    $$\kappa(X/Z,K_X+B(\bm{v}))=\kappa_{\iota}(X/Z,K_X+B(\bm{v}_0))\geq 0.$$
\end{enumerate}
\end{thm}

\medskip

\noindent\textbf{Structure of the paper}. In Section \ref{sec2}, we recall some notation and tools which will be used in this paper. In Section \ref{sec3}, we show the existence of uniform pseudo-effective rational polytopes, i.e., Corollary \ref{thm: Uniform pe polytope}. In Section \ref{sec4}, we prove Theorem \ref{thm: Uniform kod polytope assume gmm intro}. In Section \ref{sec5}, we prove Theorem \ref{thm:dccdecomp}. 

\medskip

\noindent\textbf{Acknowledgement}.
The authors would like to thank Christopher D. Hacon, Junpeng Jiao, Yujie Luo, Fanjun Meng, Lingyao Xie, and Qingyuan Xue for valuable discussions and suggestions. The first named author was supported by the China post-doctoral grants BX2021269 and 2021M702925. The second named author was supported by a grant from the Simons Foundation (Grant Number 814268, MSRI) and Start-up Grant No. JIH1414011Y of Fudan University.

\section{Preliminaries}\label{sec2}
In this paper, varieties are normal quasi-projective. We may adopt the standard notation and definitions in \cite{KM98}, and will freely use them.

\subsection{Pairs and singularities}

\begin{defn}\label{defn: DCC and ACC}
Let $\Ii\subset\Rr$ be a set. We say that $\Ii$ satisfies the \emph{descending chain condition} (DCC) if any decreasing sequence $a_1\ge a_2 \ge \cdots \ge a_k \ge\cdots$ in $\Ii$ stabilizes. We say that $\Ii$ satisfies the \emph{ascending chain condition} (ACC) if any increasing sequence in $\Ii$ stabilizes. 
\end{defn}

Let $\Ii\subset\Rr$ be a set, $X$ a normal variety, and $B:=\sum_{i=1}^s b_iB_i$ an $\Rr$-divisor on $X$, where $B_i$ are the irreducible components of $B$. We write $B\in\Ii$ if $b_i\in\Ii$ for every $i$, and $B\geq 0$ if  $B\in [0,+\infty)$. We define $||B||:=\max_{1\leq i\leq s}\{|b_i|\}$. Let $E$ be a prime divisor on $X$. We define $\mult_EB$ to be the \emph{multiplicity} of $B$ along $E$. 

\begin{defn}\label{defn: positivity}
A \emph{sub-pair} $(X/Z,B)$ consists of varieties $X,Z$, a projective morphism $X\to Z$, and an $\Rr$-divisor $B$ such that $K_X+B$ is $\Rr$-Cartier. If $B\geq 0$, then $(X/Z,B)$ is called a \emph{pair}. If $B\in [0,1]$, then $B$ is called a \emph{boundary} on $X$. A pair is also called an \emph{$\Rr$-pair}. If $B\in\Qq$, $(X/Z,B)$ is called a \emph{$\Qq$-pair}. If either $\dim Z=0$ or $Z$ is clear from the context, then we may omit $Z$.
	
Let $\phi:W\to X$ be any log resolution of $(X,B)$, and let
$$K_W+B_W:=\phi^{*}(K_X+B).$$
The \emph{log discrepancy} of a prime divisor $D$ on $W$ with respect to $(X,B)$ is $1-\mult_{D}B_W$ and is denoted by $a(D,X,B).$ For any non-negative real number $\epsilon$, we say that $(X,B)$ is \emph{lc} (resp. \emph{klt}, \emph{$\epsilon$-lc}, \emph{$\epsilon$-klt}) if $a(D,X,B)\geq 0$, (resp. $>0$, $\geq\epsilon$, $>\epsilon$) for every prime divisor $D$ over $X$. We say that $(X,B)$ is \emph{plt} (resp. \emph{$\epsilon$-plt}), if $a(D,X,B)>0$, (resp. $>\epsilon$) for every prime divisor $D$ which is exceptional over $X$.
	
We say that $(X,B)$ is \emph{dlt} if $a(D,X,B)>0$ for any exceptional prime divisor $D\subset W$ over $X$ for some log resolution $\phi:W\rightarrow X$.
\end{defn}

\begin{defn}
Let $X\rightarrow Z$ be a projective morphism between varieties and $D$ an $\Rr$-divisor on $X$. We define
$$|D/Z|:=\{D'\mid 0\leq D'\sim_Z D\},\text{ and }|D/Z|_{\Rr}:=\{D'\mid 0\leq D'\sim_{\Rr,Z} D\}.$$
If $|D/Z|_{\Rr}\not=\emptyset$, we say that $D$ is \emph{effective} over $Z$. If $Z=\{pt\}$, we may omit $Z$.
\end{defn}

\begin{defn}[Lc thresholds] 
Assume that $(X,B)$ is an lc pair. The \emph{lc threshold} of an $\Rr$-Cartier $\Rr$-divisor $G\geq 0$ with respect to $(X,B)$ is 
	$$\lct(X,B;G):=\sup\{c\ge 0\mid (X,B+cG) \text{ is lc}\}.$$ 
\end{defn}

\subsection{Iitaka dimensions and invariant Iitaka dimensions}
\begin{defn}[Iitaka dimensions, cf. {\cite[II 3.2 Definition]{Nak04}}]\label{defn: Iitaka dimension}
Let $X$ be a normal projective variety and $D$ an $\Rr$-divisor on $X$. For any positive integer $m$ such that $|\lfloor mD\rfloor|\not=\emptyset$, we define $$\Phi_m: X\dashrightarrow\Pp(H^0(X,\lfloor mD\rfloor)).$$
The \emph{Iitaka dimension} $\kappa(X,D)$ of $D$ is defined in the following way. If $|\lfloor mD\rfloor|\not=\emptyset$ for some positive integer $m$, then
$$\kappa(X,D):=\max\{\dim\Phi_m(X)\mid {m\in\Zz^+,|\lfloor mD\rfloor|\not=\emptyset}\}.$$
Otherwise, let $\kappa(X,D):=-\infty$. Note that if $|\lfloor mD\rfloor|\not=\emptyset$, then by \cite[II 3.8 Corollary]{Nak04},
$$\kappa(X,D)=\max\left\{k\in\mathbb N^+\mid \underset{m\rightarrow+\infty}{\lim\sup}\frac{\dim H^0(X,\lfloor mD\rfloor)}{k^m}>0\right\}.$$
\end{defn}

\begin{defn}[Relative Iitaka dimensions, cf. {\cite[Definition 2.1.1]{Cho08}}]
Let $f: X\rightarrow Z$ be a projective morphism between varieties, $D$ an $\Rr$-divisor on $X$, and $F$ a very general fiber of the Stein factorialization of $f$. We define the \emph{relative Iitaka dimension} $\kappa(X/Z,D)$ of $D$ over $Z$ in the following way.
\begin{enumerate}
    \item $\kappa(X/Z,D):=\kappa(F,D|_F)$, if $|\lfloor mD\rfloor/Z|\not=\emptyset$ for some positive integer $m$ and $\dim F>0$.
    \item $\kappa(X/Z,D):=-\infty$, if $|\lfloor mD\rfloor/Z|=\emptyset$ for every positive integer $m$.
    \item $\kappa(X/Z,D):=0$, if $\dim F=0$.
\end{enumerate}
\end{defn}
\begin{rem}
Most references we can find only have parts (1) and (2) as the definition. We add part (3) to make the definition rigorous. Notice that when $\dim X=\dim Z$, $|\lfloor mD\rfloor/Z|\not=\emptyset$ for any positive integer $m$.
\end{rem}

\begin{defn}[Invariant Iitaka dimensions, cf. {\cite[Definition 2.2.1]{Cho08}}]
Let $f: X\rightarrow Z$ be a projective morphism between varieties and $D$ an $\Rr$-divisor on $X$. The \emph{(relative) invariant Iitaka dimension} $\kappa_{\iota}(X/Z,D)$ of $D$ over $Z$ is defined as follows. If $|D/Z|_{\Rr}\not=\emptyset$, then we define
$$\kappa_{\iota}(X/Z,D):=\kappa(X/Z,D')$$
for some $\Rr$-divisor $D'\in|D/Z|_\Rr$. Otherwise let $\kappa_{\iota}(X/Z,D):=-\infty$. If $Z=\{pt\}$, then we may omit $Z$. Note that $\kappa_{\iota}(X/Z,D)$ is independent of the choice of $D'$ \cite[Corollary 2.1.4]{Cho08}.
\end{defn}

\begin{prop}[{\cite[Proposition 2.2.2]{Cho08}}]\label{prop: basic properties of three iitaka dimensions}
Let $X\rightarrow Z$ be a projective morphism between varieties, and $D$ and $D'$ two $\Rr$-divisors on $X$. Then we have the following:
\begin{enumerate}
    \item $\kappa(X/Z,D)\leq\kappa_{\iota}(X/Z,D)$, and $\kappa(X/Z,D)<\kappa_{\iota}(X/Z,D)$ if and only if $\kappa(X/Z,D)=-\infty$ and $\kappa_{\iota}(X/Z,D)\geq 0$.
    \item If $D\sim_{\Rr,Z}D'$, then $\kappa_{\iota}(X/Z,D)=\kappa_{\iota}(X/Z,D')$.
\end{enumerate}
\end{prop}

The following example shows that we may have strict inequality in Proposition \ref{prop: basic properties of three iitaka dimensions}(1).

\begin{ex}\label{ex: why invairant iitaka dimension}
Let $X:=\Pp^1$ and $p_1,p_2,p_3,p_4$ four different closed points on $X$. For any real number $a$, we let $B(a):=a(p_1+p_2)+(1-a)(p_3+p_4)$. Then whenever $a\in (0,1)$, $(X,B(a))$ is klt and $K_X+B(a)\sim_{\Rr}0$. However, for any integer $m$,
$$\deg \lfloor m(K_X+B(a))\rfloor=-2,\text{ if }ma\not\in\Zz,$$
and
$$\deg \lfloor m(K_X+B(a))\rfloor=0,\text{ if }ma\in\Zz.$$
This implies that
$\kappa(X,K_X+B(a))=-\infty$ when $a\not\in\Qq$ and $\kappa(X,K_X+B(a))=0$ when $a\in\Qq$. Nevertheless, $\kappa_{\iota}(X,K_X+B(a))=0$ for any $a\in\Rr$.
\end{ex}

\begin{lem}\label{lem: kappaiota is kappa for q divisors}
Let $X\rightarrow Z$ be a projective morphism between varieties and $D$ a $\Qq$-divisor on $X$. Then $\kappa_{\iota}(X/Z,D)=\kappa(X/Z,D)$.
\end{lem}
\begin{proof}
We may assume that $\dim X>\dim Z$, otherwise there is nothing to prove. Suppose on the contrary that $\kappa_{\iota}(X/Z,D)\not=\kappa(X/Z,D)$. Then $\kappa_{\iota}(X/Z,D)\geq 0$ and $\kappa(X/Z,D)=-\infty$ by Proposition \ref{prop: basic properties of three iitaka dimensions}(1). Since $\kappa_{\iota}(X/Z,D)\geq 0$, $D\sim_{\Rr,Z}D'\geq 0$ for some $\Rr$-divisor $D'$ on $X$.

We may write $D':=D_0+\sum_{i=1}^cr_iD_i$ for some $\Qq$-divisors $D_0,\dots,D_c$, and real numbers $r_1,\dots,r_c$ such that $1,r_1,\dots,r_c$ are linearly independent over $\Qq$. Moreover, as $D'\geq 0$, we may pick rational numbers $r_1',\dots,r_c'$ that are sufficiently close to $r_1,\dots,r_c$ respectively, such that $$D'':=D_0+\sum_{i=1}^cr_i'D_i\geq 0.$$
Since $D$ is a $\Qq$-divisor and $D\sim_{\Rr}D'$, we see that $D\sim_{\Qq,Z}D_0$ and $D_i\sim_{\Qq,Z}0$ for every $1\le i\le c$ by \cite[Lemma 5.3]{HLS19}. It follows that $D\sim_{\Qq,Z}D''\geq 0$. Hence one can find a positive integer $m$ such that $mD,mD''$ are integral and $mD\sim_Z mD''$. Thus $\kappa(X/Z,D)\geq 0$, a contradiction.
\end{proof}

\subsection{Pseudo-effective thresholds}

\begin{defn}[Pseudo-effective thresholds]
Let $(X/Z,B)$ be a pair, and $D\geq 0$ an $\Rr$-Cartier $\Rr$-divisor on $X$. We define
\begin{align*}
\pet(X/Z,B;D):=\inf\left\{\begin{array}{l|l}
+\infty,t& \begin{array}{l}
t\geq 0,\  (X,B+tD)\text{ is lc, and}\\ K_X+B+tD\text{ is pseudo-effective over }Z
\end{array}
\end{array}\right\}
\end{align*}
to be the \emph{pseudo-effective threshold} of $D$ with respect to $(X/Z,B)$. 

For any positive integer $d$, and sets $\Ii\subset [0,1]$ and $\Ii'\subset [0,+\infty)$, we define
\begin{align*}
\PET(d,\Ii;\Ii'):=\left\{\begin{array}{l|l}
\pet(X/Z,B;D)& \begin{array}{l}
(X/Z,B)\text{ is lc, }\dim X=d,\ B\in\Ii,\text{ and}\\
     D\in\Ii'\text{ is an }\Rr\text{-Cartier }\Rr\text{-divisor on }X
\end{array}
\end{array}\right\}.
\end{align*}
\end{defn}

\begin{thm}[ACC for PETs, cf. {\cite[Theorem 1.6]{DC16}}]\label{thm: acc pet}
Let $d$ be a positive integer, and $\Ii\subset [0,1]$ and $\Ii'\subset [0,+\infty)$ two DCC sets. Then $\PET(d,\Ii;\Ii')$ satisfies the ACC.
\end{thm}
The proof essentially follows from \cite[proof of Theorem 1.6]{DC16} and \cite[Theorem 1.4]{HMX14} when $Z$ is a point. For the reader's convenience, we give a full proof here. 
\begin{proof}
Possibly replacing $\Ii$ with $\Ii\cup\{1\}$ we may assume that $1\in\Ii$. Suppose that the theorem does not hold. Then there exist a sequence of lc pairs $(X_i/Z_i,B_i)$ of dimension $d$, and a sequence of $\Rr$-Cartier $\Rr$-divisors $D_i$ on $X$ such that $B_i\in\Ii$, $D_i\in\Ii'$, and $t_i:=\pet(X_i/Z_i,B_i;D_i)$ is strictly increasing. According to \cite[Theorem 1.1]{HMX14}, possibly passing to a subsequence, we may assume that $0<t_i<\lct(X_i,B_i;D_i)$ for every $i$.

Let $f_i: W_i\rightarrow X_i$ be a dlt modification of $(X_i,B_i)$,
$$K_{W_i}+B_{W_i}:=f_i^*(K_{X_i}+B_i),\text{ and }D_{W_i}:=f_i^*D_i.$$
Then $B_{W_i}\in\Ii$ and $f_i^*D_i=(f_i^{-1})_*D_i$. It is clear that $t_i=\pet(W_i/Z_i,B_{W_i};D_{W_i})$. Possibly replacing $X_i,B_i$ and $D_i$ with $W_i,B_{W_i}$ and $D_{W_i}$ respectively, we may assume that each $X_i$ is $\Qq$-factorial.

Let $t_0:=0$. Since $t_i$ is strictly increasing, for every positive integer $i$, we may pick a real number $t_{i-1}<h_i<t_i$, such that $\lim_{i\rightarrow+\infty}(t_i-h_i)=0$. By \cite[Lemma 5.21]{HLS19}, we may assume that $h_i$ is strictly increasing. Now $(X_i,B_i+h_iD_i)$ is lc, and $K_{X_i}+B_i+h_iD_i$ is not pseudo-effective over $Z_i$, one can run a $(K_{X_i}+B_i+h_iD_i)$-MMP with scaling of some ample divisor over $Z_i$ which terminates with a Mori fiber space $f_i: Y_i\rightarrow T_i$ over $Z_i$. Let $B_{Y_i}$ and $D_{Y_i}$ be the strict transforms of of $B_i$ and $D_i$ on $Y_i$ respectively, and $s_i:=\lct(Y_i,B_{Y_i};D_{Y_i})$. 
\begin{claim}\label{claim: lct larger than original}
$s_i\geq t_i$ except for finitely many $i$.
\end{claim}
\begin{proof}[Proof of Claim \ref{claim: lct larger than original}]
Suppose that the claim does not hold. Then possibly passing to a subsequence, we may assume that $s_i<t_i$ for every $i$. Since $(Y_i,B_{Y_i}+h_iD_{Y_i})$ is lc, $s_i\geq h_i$. Since $\lim_{i\rightarrow+\infty}(t_i-h_i)=0$, $\lim_{i\rightarrow+\infty}(t_i-s_i)=0$. By \cite[Lemma 5.21]{HLS19}, there is a subsequece of $s_i$ that is strictly increasing, which contradicts \cite[Theorem 1.4]{HMX14}. Thus the claim holds.
\end{proof}
\noindent\textit{Proof of Theorem \ref{thm: acc pet} continued.} 
By Claim \ref{claim: lct larger than original}, possibly passing to a subsequence, we may assume that $s_i\geq t_i$ for every $i$. In particular, $(Y_i,B_{Y_i}+t_iD_{Y_i})$ is lc for every $i$. Since $K_{X_i}+B_i+t_iD_i$ is pseudo-effective over $Z_i$, $K_{Y_i}+B_{Y_i}+t_iD_{Y_i}$ is pseudo-effective over $Z_i$ and thus $K_{Y_i}+B_{Y_i}+t_iD_{Y_i}$ is nef over $T_i$. Moreover, as $K_{Y_i}+B_{Y_i}+h_iD_{Y_i}$ is anti-ample over $T_i$, $D_{Y_i}\not\equiv_{T_i}0$, and there exists a real number $\eta_i\in (h_i,t_i]$, such that $K_{Y_i}+B_{Y_i}+\eta_iD_{Y_i}\equiv_{T_i}0$. It is clear that $(Y_i,B_{Y_i}+\eta_iD_{Y_i})$ is lc and $\eta_i$ is strictly increasing. Let $F_i$ be a general fiber of $f_i$, and
$$K_{F_i}+\Delta_i:=(K_{Y_i}+B_{Y_i}+\eta_iD_{Y_i})|_{F_i}.$$
Then $K_{F_i}+\Delta_i\equiv 0$, $(F_i,\Delta_i)$ is lc, and the set of coefficients of $\Delta_i$ is a DCC but not finite set, which contradicts \cite[Theorem 1.4]{HMX14}. This completes the proof.
\end{proof}

Then next corollary is a variation of Theorem \ref{thm: acc pet} which is more convenient to apply in many situations.
\begin{cor}\label{cor: gap of psd thresholds}
Let $d$ be a positive integer and $\Ii\subset [0,1]$ a DCC set. Then there exist a finite set $\Ii_0\subset\bar\Ii$ and a function $h: \bar\Ii\rightarrow\Ii_0$ depending only on $d$ and $\Ii$ satisfying the following. 
\begin{enumerate}
    \item For any $\gamma\in\bar\Ii$, $0\leq h(\gamma)\leq\gamma$, and if $\gamma>0$, then $h(\gamma)>0$. 
    \item Assume that $(X/Z,B:=\sum b_iB_i\in\Ii)$ is an lc pair of dimension $d$ such that $B\in\Ii$ and $B_i$ are the irreducible components of $B$. Then 
\begin{enumerate}
  \item if $K_X+B$ is pseudo-effective (resp. big) over $Z$,  then $K_X+\sum h(b_i)B_i$ is pseudo-effective (resp. big) over $Z$, and 
  \item if $B$ is big over $Z$, then $\sum h(b_i)B_i$ is big over $Z$.
\end{enumerate}
\end{enumerate}
\end{cor}
\begin{proof}
Possibly replacing $\Ii$ with $\Ii\cup\{1\}$ and $(X,B)$ with a dlt modification, we may assume that $X$ is $\Qq$-factorial. According to Theorem \ref{thm: acc pet}, the set
$$\Ii':=\PET(d,\bar\Ii;\bar\Ii)$$
satisfies the ACC. Then we can define
$$\beta:=\max\{\frac{1}{2},\frac{\gamma+1}{2}\mid\gamma\in\Ii',\gamma<1\}\in (0,1).$$
Since $\Ii$ satisfies the DCC, we may define
$$\gamma_0:=\min\{1,\gamma\mid\gamma>0,\gamma\in\Ii\}>0\text{ and }\epsilon:=\frac{(1-\beta)\gamma_0}{\beta}.$$
By \cite[Theorem 5.18]{HLS19}, there exist a finite set $\Ii_0'\subset\bar\Ii$  and a projection $g:\bar\Ii\rightarrow\Ii_0'$ depending only on $d$ and $\Ii$, such that $\gamma+\epsilon\geq g(\gamma)\geq\gamma$ for any $\gamma\in\bar\Ii$, and $(X,\sum g(b_i)B_i)$ is lc. Now we may define $h$ as follows. For any $\gamma\in\bar\Ii$, if $g(\gamma)=\gamma$, then $h(\gamma):=\gamma$, otherwise let $h(\gamma):=\beta g(\gamma)$.

We show that $\Ii_0:=\Ii_0'\cup\beta\Ii_0'$ and $h$ satisfy our requirements. It is clear that $h(\gamma)\geq 0$ for any $\gamma\in\bar\Ii$ and $h(\gamma)>0$ when $\gamma>0$. For any $\gamma\in\bar\Ii$, if $g(\gamma)=\gamma$, then $h(\gamma)=\gamma\leq\gamma$ by construction, and if $g(\gamma)>\gamma$, then 
$$h(\gamma)=\beta g(\gamma)\leq \beta(\gamma+\epsilon)=\beta\gamma+(1-\beta)\gamma_0\leq\beta\gamma+(1-\beta)\gamma=\gamma,$$
which implies (1). 

Let
$$F:=\sum_{g(b_i)=b_i}g(b_i)B_i \text{ and } G:=\sum_{g(b_i)>b_i}g(b_i)B_i.$$ 
If $K_X+B$ is pseudo-effective (resp. big) over $Z$, then by our construction, $\pet(X/Z,F;G)\in\PET(d,\bar\Ii;\bar\Ii)$ and $\pet(X/Z,F;G)<1$. Hence 
$$\pet(X/Z,F;G)\leq 2\beta-1<\beta.$$ 
Since
\begin{align*}
    K_X+\sum h(b_i)B_i
    =K_X+F+\beta G
    =\frac{1}{2}(K_X+F+(2\beta-1)G)+\frac{1}{2}(K_X+F+G),
\end{align*}
we see that $K_X+\sum h(b_i)B_i$ is pseudo-effective (resp. big) over $Z$. If $B$ is big over $Z$, then $\sum h(b_i)B$ is big over $Z$ as $\Supp B=\sum\Supp h(b_i)B$. We may finish the proof.
\end{proof}

\section{Uniform pseudo-effective rational polytopes}\label{sec3}
In this section, we show the existence of uniform pseudo-effective rational polytopes.

\begin{thm}\label{thm: one variable psd polytope}
Let $d$ be a positive integer, $r_1,\dots,r_c$ real numbers, and $s_1,\dots,s_m:\Rr^{c+1}\rightarrow\Rr$ $\Qq$-linear functions such that $1,r_1,\dots,r_c$ are linearly independent over $\Qq$. Let $\bm{r}(t):=(r_1,\dots,r_{c-1},t)\in\Rr^c$ for any real number $t$. Then there exists a positive real number $\epsilon$ depending only on $d,r_1,\dots,r_c$ and $s_1,\dots,s_m$ satisfying the following. 

Assume that $X$ is a variety of dimension $d$, $B_1,\dots,B_m\geq 0$ are Weil divisors on $X$, and $X\to Z$ is a projective morphism between varieties. Let $B(t):=\sum_{j=1}^ms_j(1,\bm{r}(t))B_j$ for any real number $t$. Suppose that $(X,B(r_c))$ is lc and $K_X+B(r_c)$ is pseudo-effective over $Z$. Then $(X,B(t))$ is lc and $K_X+B(t)$ is pseudo-effective over $Z$ for any $t$ such that $|t-r_c|<\epsilon$.
\end{thm}
\begin{proof}
First we reduce to the case when $X$ is $\Qq$-factorial. Let $h: W\rightarrow X$ be a dlt modification of $(X,B(t_c))$. Let $B_W(t):=\sum_{j=1}^ms_j(1,\bm{r}(t))B_{j,W}+E$ for any real number $t$, where $E$ is the reduced exceptional divisor of $h$ and $B_{j,W}$ is the strict transform of $B_j$ on $W$ for any $1\le j\le m$. Possibly replacing $s_1,\dots,s_{m-1},s_m$ with $s_1,\dots,s_{m},1$, and $B_1,\dots,B_{m-1},B_m$ with $B_{1,W},\dots,B_{m,W},E$, we may assume that $X$ is $\Qq$-factorial.

Assume that the theorem does not hold. Then there exist a sequence of $\Qq$-factorial varieties $X^i$ of dimension $d$, Weil divisors $B_1^i,\dots,B_m^i\ge0$ on $X^i$, projective morphisms $X^i\to Z^i$ between varieties, and real numbers $t_i$ such that $(X^i/Z^i,B^i(r_c))$ is lc, $K_{X^i}+B^i(r_c)$ is pseudo-effective over $Z^i$, $K_{X^i}+B^i(t_i)$ is not pseudo-effective over $Z^i$, and $|t_i-r_c|<\frac{1}{i}$, where $B^i(t):=\sum_{j=1}^ms_j(1,\bm{r}(t))B^i_j$ for each $i$ and for any real number $t$. According to \cite[Corollary 5.5]{HLS19}, possibly passing to a subsequence, we may assume that  $(X^i,B^i(t_i))$ is lc for each $i$. Without lost of generality, we may assume that $r_c>t_i>t_{i-1}$ for any $i$.

We may run a $(K_{X^i}+B^i(t_i))$-MMP with scaling of some ample divisor over $Z^i$, which terminates with a Mori fiber space $f_i: Y^i\rightarrow T^i$ over $Z^i$. For every $i,j$ and any real number $t$, let $B^i_{j,Y}$ and $B^i_Y(t)$ be the strict transforms of $B^i_j$ and $B^i(t)$ on $Y^i$ respectively. We may write
$$B^i_Y(t)=B^i_Y(r_c)-(r_c-t)F_i+(r_c-t)G_i$$
for some $\Qq$-divisors $F_i\geq 0$ and $G_i\geq 0$ such that $F_i\wedge G_i=0$. Since $(Y^i,B^i_Y(t_i))$ is lc, $(Y^i,B^i_Y(r_c)-(r_c-t_i)F_i)$ is lc. Since $\lim_{i\rightarrow+\infty}t_i=r_c$, $r_c>t_i$ for each $i$, and $B^i_Y(r_c)\in\{s_j(r_c)\mid 1\le j\le m\}$ which is a finite set, by \cite[Theorem 1.4]{HMX14}, possibly passing to a subsequence, we may assume that $(Y^i,B^i_Y(r_c))$ is lc.

Since $K_{X^i}+B^i(r_c)$ is pseudo-effective over $Z^i$ and $K_{Y^i}+B^i_Y(t_i)$ is anti-ample over $T^i$, there exists a real number $\eta_i\in (t_i,r_c]$ such that $K_{Y^i}+B^i_Y(t_i)\equiv_{T^i}0$. If $\eta_i=r_c$, by \cite[Lemma 5.3]{HLS19}, $K_{Y^i}+B^i_Y(t)\equiv_{T^i}0$ for any real number $t$, a contradiction. Hence $\eta_i<r_c$. 
Let $F^i$ be a general fiber of $f_i$, $B^i_{j,F}:=B^i_{j,Y}|_{F^i}$, and $B^i_F(t):=\sum_{j=1}^ms_j(1,\bm{r}(t))B^i_{j,F}$ for each $i$ and for any real number $t$. Then $(F^i,B^i_{F}(\eta_i))$ is lc,
$$K_{F^i}+B^i_{F}(\eta_i)=(K_{Y^i}+B^i_Y(\eta_i))|_{F^i}\equiv 0,$$
and $K_{F^i}+B^i_{F}(t_i)$ is anti-ample. However, this contradicts \cite[Theorem 3.8]{Nak16} as $\eta_i<r_c$ and $\lim_{i\rightarrow+\infty}\eta_i=r_c$.
\end{proof}

\begin{cor}\label{thm: Uniform pe polytope}
Let $d$ be a positive integer, $\bm{v}_0:=(v^0_1,\ldots,v^0_m)\in\Rr^m$ a point, and $V\subset\Rr^m$ the rational envelope of $\bm{v}_0$. Then there exists an open set $V_0\ni\bm{v}_0$ of $V$ depending only on $d$ and $\bm{v}_0$ satisfying the following. 

Assume that $X$ is a variety of dimension $d$, $B_1,\dots,B_m\geq 0$ are Weil divisors on $X$, and $X\to Z$ is a projective moprhism between varieties. Let $B(\bm{v}):=\sum_{i=1}^mv_iB_m$ for any point $\bm{v}:=(v_1,\dots,v_m)\in\Rr^m$. Suppose that $(X,B(\bm{v}_0))$ is lc and $K_X+B(\bm{v}_0)$ is pseudo-effective over $Z$. Then $(X,B(\bm{v}))$ is lc and $K_X+B(\bm{v})$ is pseudo-effective over $Z$ for any $\bm{v}\in V_0$.
\end{cor}

\begin{proof}
There exist real numbers $r_1,\ldots,r_c$ depending only on $d$ and $\bm{v}_0$, such that $\{1,r_1,\ldots,r_{c}\}$ is a basis of $\Span_{\Qq}(\{1,v^0_1,\ldots,v^0_{m}\})$ over $\Qq$. When $c=0$, we may let $V_0:=V=\{\bm{v}\}$, and the theorem holds. In the following, we may assume that $c\geq 1$. 

Let $\bm{r}(t):=(r_1,\dots,r_{c-1},t)$ for any real number $t$. There exist $\Qq$-linear functions $s_i:\Rr^{c+1}\to \Rr$, such that $s_i(1,\bm{r}(r_c))=v_i$ for any $1\le i\le m$, and the map
$$(s_1(1,\bm{x}),\ldots,s_m(1,\bm{x})):\Rr^{c}\to V$$
is one-to-one. It suffices to show that there exists an open set $U_c\ni\bm{r}(r_c)$ of $\Rr^{c}$, such that for any $\bm{x}\in U_c$, $(X,\sum_{i=1}^m s_i(1,\bm{x})B_{i})$ is lc, and $K_X+\sum_{i=1}^m s_i(1,\bm{x})B_{i}$ is pseudo-effective over $Z$.

We prove Theorem \ref{thm: Uniform pe polytope} by induction on $c$. When $c=1$, the theorem follows from Theorem \ref{thm: one variable psd polytope}. Assume that $c\ge 2$. For simplicity, for every $1\leq i\leq m$, we let
$$s_i(t):=s_i(1,\bm{r}(t)).$$ 
According to Theorem \ref{thm: one variable psd polytope}, there 
exist two positive real numbers $\epsilon_1$ and $\epsilon_2$ depending only on $d$, $r_i$, $s_j$, such that $r_c+\epsilon_1,r_c-\epsilon_2\in\Qq$,
\begin{itemize}
    \item both $(X,\sum_{i=1}^m s_i(r_c+\epsilon_1)B_i)$ and $(X,\sum_{i=1}^m s_i(r_c-\epsilon_2)B_i)$ are lc, and
    \item both $K_X+\sum_{i=1}^m s_i(r_c+\epsilon_1)B_i$ and $K_X+\sum_{i=1}^m s_i(r_c-\epsilon_2)B_i$ are pseudo-effective over $Z$.
\end{itemize}
By induction, there exists an open set $U_{c-1}\ni (r_1,\ldots,r_{c-1})$ of $\Rr^{c-1}$, such that
\begin{itemize}
    \item both $(X,\sum_{i=1}^m s_i(1,\bm{x}',r_{c}+\epsilon_1)B_{i})\text{ and }(X,\sum_{i=1}^m s_i(1,\bm{x}',r_{c}-\epsilon_2)B_{i})$ are lc, and
    \item both $K_X+\sum_{i=1}^m s_i(1,\bm{x}',r_{c}+\epsilon_1)B_{i}\text{ and }K_X+\sum_{i=1}^m s_i(1,\bm{x}',r_{c}-\epsilon_2)B_{i}$ are pseudo-effective over $Z$,
\end{itemize}
for any $\bm{x}':=(x_1,\ldots,x_{c-1})\in U_{c-1}$. We get the desired $U_c$ by letting $U_c:=U_{c-1}\times (r_{c}-\epsilon_2,r_{c}+\epsilon_1)$.
\end{proof}

\section{Uniform rational polytopes for Iitaka dimensions}\label{sec4}
\subsection{Iitaka dimensions and perturbation theory}

\begin{lem}\label{lem: real kodaira dimension positive has support irrational}
Assume that $X\rightarrow Z$ is a projective morphism between normal quasi-projective varieties, and $D:=D_0+D_1$ is an $\Rr$-divisor on $X$ such that $D_0\in\Qq$, $D_1\in\Rr\setminus\Qq$ and $\kappa(X/Z,D)\geq 0$. Then there exists a positive integer $m_0$, such that for any positive integer $m$ divisible by $m_0$, there exists a Weil divisor $L_m\in |\lfloor mD\rfloor/Z|$ such that $\Supp D_1\subset\Supp L_m$.
\end{lem}
\begin{proof}
By assumption, we may find a positive integer $m_1$ such that $m_1D_0\in\Zz$ and $|\lfloor m_1D\rfloor/Z|\not=\emptyset$. It is clear that $\Supp D_1=\Supp \{m_1D_1\}$, and
$$\delta:=\min\{\mult_C\{m_1D_1\}\mid C\subset\Supp D_1\text{ is a prime divisor}\}>0.$$
Let $m_0:=\lceil\frac{1}{\delta}\rceil m_1$, then
$$\lfloor m_0D\rfloor=\lfloor\lceil\frac{1}{\delta}\rceil m_1D\rfloor=\lfloor\lceil\frac{1}{\delta}\rceil(\lfloor m_1D\rfloor +\{m_1D\})\rfloor=\lceil\frac{1}{\delta}\rceil\lfloor m_1D\rfloor+\lfloor\lceil\frac{1}{\delta}\rceil\{m_1D\}\rfloor.$$
We may pick $0\leq P\in |\lfloor m_1D\rfloor/Z|$ and set $L_{m_0}:=\lceil\frac{1}{\delta}\rceil P+\lfloor\lceil\frac{1}{\delta}\rceil\{m_1D\}\rfloor.$ It is clear that $L_{m_0}\in|\lf m_0D\rf/Z|$ and $\Supp D_1\subset\Supp L_{m_0}.$ For any positive integer $m$ divisible by $m_0$, we may let
$$L_m:=\frac{m}{m_0}L_{m_0}+(\lfloor mD\rfloor-\frac{m}{m_0}\lfloor m_0D\rfloor).$$
Then $0\le L_m\sim_Z\lfloor mD\rfloor$ and $\Supp D_1\subset\Supp L_m$. This finishes the proof.
\end{proof}

\begin{lem}\label{lem: rational kod geq 0 and central kod geq 0 imply kod geq 0 on open set}
Let $\bm{r}:=(r_1,\dots,r_c)\in\Rr^c$ be a point, and $U\ni\bm{r}$ an open subset of $\Rr^c$ such that $1,r_1,\dots,r_c$ are linearly independent over $\Qq$. Assume that $X\rightarrow Z$ is a projective morphism between varieties, and $D_0,D_1,\dots,D_c$ are $\Qq$-divisors on $X$. Let $D(\bm{v}):=D_0+\sum_{i=1}^cv_iD_i$ for any $\bm{v}:=(v_1,\dots,v_c)\in\Rr^c$. Suppose that $\kappa(X/Z,D(\bm{r}))\geq 0$, and $\kappa(X/Z,D(\bm{v}))\geq 0$ for any point $\bm{v}\in U\cap\Qq^c$. Then $\kappa(X/Z,D(\bm{v}))\geq 0$ for any point $\bm{v}\in U$.
\end{lem}
\begin{proof}
We may assume that $\dim X>\dim Z$, otherwise there is nothing to prove. Fix a point $\bm{v}:=(v_1,\dots,v_c)\in U\setminus\Qq^c$. By Lemma \ref{lem: real kodaira dimension positive has support irrational}, there exist a positive integer $m_1$ and a Weil divisor $L_1\geq 0$ on $X$, such that $L_1\in |\lfloor m_1D(\bm{r})\rfloor/Z|$, $m_1D_0$ is a Weil divisor, and $\cup_{i=1}^c\Supp D_i\subset\Supp L_1$. Since $U$ is an open set, there exists an integer $m_0\geq 2$ such that $$\bm{u}:=\bm{r}+\frac{m_0}{m_0-1}(\bm{v}-\bm{r})\in U.$$
We may pick a point $\bm{u}_0\in U\cap\Qq^c$, such that $$||D(\bm{u})-D(\bm{u}_0)||<\frac{1}{m_1(m_0-1)}.$$
By construction,
$$D(\bm{v})=\frac{1}{m_0}D(\bm{r})+\frac{m_0-1}{m_0}D(\bm{u}),$$
which implies that
 $$||D(\bm{v})-(\frac{1}{m_0}D(\bm{r})+\frac{m_0-1}{m_0}D(\bm{u}_0))||<\frac{1}{m_1m_0}.$$
Equivalently,
$$||m_1m_0D(\bm{v})-(m_1D(\bm{r})+m_1(m_0-1)D(\bm{u}_0))||<1.$$
Since $\bm{u}_0\in U\cap\Qq^c$, $\kappa(X/Z,D(\bm{u}_0))\geq 0$. Thus we may pick a positive integer $m_2$, such that $m_2D(\bm{u}_0)$ is an integral divisor and $|m_2D(\bm{u}_0)/Z|\not=\emptyset$. We may pick $L_2\in |m_2D(\bm{u}_0)/Z|$. 

Let $E:=\sum_{i=1}^c\Supp D_i$. By our construction, $L_1\geq E$. Since
$$\Supp\left(m_1m_0D(\bm{v})-(m_1D(\bm{r})+m_1(m_0-1)D(\bm{u}_0))\right)\subset\cup_{i=1}^c\Supp D_i$$
and
$||m_1m_0D(\bm{v})-(m_1D(\bm{r})+m_1(m_0-1)D(\bm{u}_0))||<1,$
we have
$$m_1m_0D(\bm{v})-(m_1D(\bm{r})+m_1(m_0-1)D(\bm{u}_0))\geq -E.$$
It follows that
$$\lfloor m_2m_1m_0D(\bm{v})-m_2(m_1D(\bm{r})+m_1(m_0-1)D(\bm{u}_0))\rfloor \geq -m_2E.$$
Moreover, we have
\begin{align*}
    &\ \ \ \lfloor m_2m_1D(\bm{r})+m_2m_1(m_0-1)D(\bm{u}_0)-m_2E\rfloor=m_2m_1(m_0-1)D(\bm{u}_0)+\lfloor m_2(m_1D(\bm{r})-E)\rfloor\\&
    \geq m_2m_1(m_0-1)D(\bm{u}_0)+m_2\lfloor(m_1D(\bm{r})-E)\rfloor\sim_Z m_1(m_0-1)L_2+m_2(\lfloor m_1D(\bm{r})\rfloor-E)\\&
    \sim_Z m_1(m_0-1)L_2+m_2(L_1-E)\geq 0.
\end{align*}
Thus we may pick $L_3\in |\lfloor m_2m_1D(\bm{r})+m_2m_1(m_0-1)D(\bm{u}_0)-m_2E\rfloor/Z|.$
We have
\begin{align*}
    &\lfloor m_2m_1m_0D(\bm{v})\rfloor\\
    =&\lfloor m_2m_1m_0D(\bm{v})-(m_2m_1D(\bm{r})+m_2m_1(m_0-1)D(\bm{u}_0))\\
    &+(m_2m_1D(\bm{r})+m_2m_1(m_0-1)D(\bm{u}_0))-m_2E+m_2E\rfloor\\
    \geq& \lfloor m_2m_1m_0D(\bm{v})-(m_2m_1D(\bm{r})+m_2m_1(m_0-1)D(\bm{u}_0))\rfloor\\
    &+\lfloor (m_2m_1D(\bm{r})+m_2m_1(m_0-1)D(\bm{u}_0))-m_2E\rfloor+m_2E\\
    \geq &-m_2E+\lfloor ((m_2m_1D(\bm{r})+m_2m_1(m_0-1)D(\bm{u}_0))-m_2E\rfloor+m_2E\\
    =&\lfloor (m_2m_1D(\bm{r})+m_2m_1(m_0-1)D(\bm{u}_0))-m_2E\rfloor\sim_Z L_3\geq 0.
\end{align*}
Therefore $|\lfloor m_2m_1m_0D(\bm{v})\rfloor/Z|\neq0$. In particular, $\kappa(X/Z,D(\bm{v}))\geq 0$.
\end{proof}

\begin{lem}\label{lem: inv kod linear geq maximum}
Let $a,b$ be two positive real numbers, $X\rightarrow Z$ a projective morphism between varieties, and $A,B,C$ three $\Rr$-Cartier $\Rr$-divisors on $X$ such that $C=aA+bB$, $\kappa_{\iota}(X/Z,A)\ge0$ and $\kappa_{\iota}(X/Z,B)\geq 0$. Then $$\kappa_{\iota}(X/Z,C)\geq\max\{\kappa_{\iota}(X/Z,A),\kappa_{\iota}(X/Z,B)\}.$$
\end{lem}
\begin{proof}
We may pick $0\leq A'\sim_{\Rr,Z}A$, $0\leq B'\sim_{\Rr,Z}B$, and let $C':=aA'+bB'$. Let $m$ be a positive integer such that $ma\geq 1$ and $mb\geq 1$. Then $mC'\geq A'+B'$. Thus $\kappa(X/Z,mC')\geq\kappa(X/Z,A')$ and  $\kappa(X/Z,mC')\geq\kappa(X/Z,B')$. The lemma follows from Proposition \ref{prop: basic properties of three iitaka dimensions}.
\end{proof}

\begin{lem}[{cf. \cite[Proposition 4.1]{Li22}}]\label{lem: kod is cons on open set}
Let $m$ be a positive integer, $V\subset\Rr^m$ an affine subspace, $U$ an open subset of $V$, $X\rightarrow Z$ a projective morphism between varieties, and $D_1,\dots,D_m$ $\Rr$-divisors on $X$. Let $D(\bm{v}):=\sum_{i=1}^mv_iD_i$ for any point $\bm{v}:=(v_1,\dots,v_m)\in\Rr^m$. 

Then if $\kappa_{\iota}(X/Z,D(\bm{v}))\geq 0$ (resp. $\kappa(X/Z,D(\bm{v}))\geq 0$)  for any $\bm{v}\in U$, then $\kappa_{\iota}(X/Z,D(\bm{v}))$ (resp. $\kappa(X/Z,D(\bm{v}))$) is a constant for any $\bm{v}\in U$.
\end{lem}
\begin{proof}
Let $\bm{v}_1\in U$ (resp. $\bm{v}_2\in U$) be a point such that $$\kappa_{\iota}(X/Z,D(\bm{v}_1))=\max_{\bm{v}\in U}\{\kappa_{\iota}(X/Z,D(\bm{v}))\}\text{ (resp. }\kappa(X/Z,D(\bm{v}_2))=\max_{\bm{v}\in U}\{\kappa(X/Z,D(\bm{v}))\}\text{)}.$$
Fix a point $\bm{v}\in U$. Choose a rational number $c_1>0$ (resp. $c_2>0$) such that
$$\bm{u}_1:=\bm{v}+c_1(\bm{v}-\bm{v}_1)\in U\text{ (resp. }\bm{u}_2:=\bm{v}+c_2(\bm{v}-\bm{v}_2)\in U\text{)}.$$
Then 
$$D(\bm{v})=\frac{1}{c_1+1}D(\bm{u}_1)+\frac{c_1}{c_1+1}D(\bm{v}_1)\text{ (resp. }D(\bm{v})=\frac{1}{c_2+1}D(\bm{u}_2)+\frac{c_2}{c_2+1}D(\bm{v}_2)\text{)}.$$
By Proposition \ref{prop: basic properties of three iitaka dimensions}, Lemma \ref{lem: inv kod linear geq maximum} and the choice of $\bm{v}_1$ (resp. $\bm{v}_2$),
\begin{center}
  $\kappa_{\iota}(X/Z,D(\bm{v}))\geq\kappa_{\iota}(X/Z,D(\bm{v}_1))\geq\kappa_{\iota}(X/Z,D(\bm{v}))$
\end{center}
\begin{center}
  $\big(\text{resp. } \kappa(X/Z,D(\bm{v}))=\kappa_{\iota}(X/Z,D(\bm{v}))\geq\kappa_{\iota}(X/Z,D(\bm{v}_2))=\kappa(X/Z,D(\bm{v}_2))\geq\kappa(X/Z,D(\bm{v}))\big).$
\end{center}
Thus $\kappa_{\iota}(X/Z,D(\bm{v}))=\kappa_{\iota}(X/Z,D(\bm{v}_1))$ (resp. $\kappa(X/Z,D(\bm{v}))=\kappa(X/Z,D(\bm{v}_2))$) is constant.
\end{proof}

\subsection{Proof of Theorem \ref{thm: Uniform kod polytope assume gmm intro}}

\begin{thm}\label{thm: Uniform kod polytope assume gmm intro1}
Let $d$ be a positive integer, $\epsilon>\epsilon'>0$ two real numbers, and $\bm{v}_0:=(v_1^0,\dots,v_m^0)\in\Rr^m$ a point. Then there exists an open set $V_0\ni\bm{v}_0$ of $\Rr^m$ depending only on $d,\epsilon,\epsilon'$ and $\bm{v}_0$ satisfying the following. 

Assume that $X\rightarrow Z$ is a projective morphism between normal quasi-projective varieties,  $\dim X=d$, and $B_1,\dots,B_m\geq 0$ are Weil divisors on $X$. Let $B(\bm{v}):=\sum_{i=1}^mv_iB_m$ for any point $\bm{v}:=(v_1,\dots,v_m)\in\Rr^m$. Assume that $(X,B(\bm{v}_0))$ is lc (resp. klt, $\epsilon$-lc, $\epsilon$-plt, $\epsilon$-klt), and $K_X+B(\bm{v}_0)$ is pseudo-effective over $Z$. Then for any point $\bm{v}\in V_0$,
\begin{enumerate}
     \item $(X,B(\bm{v}))$ is lc (resp. klt, $\epsilon'$-lc, $\epsilon'$-plt, $\epsilon'$-klt),
     \item $K_X+B(\bm{v})$ is pseudo-effective over $Z$, and
    \item if $K_X+B(\bm{v}_0)$ (resp. $B(\bm{v}_0)$) is big over $Z$, then $K_X+B(\bm{v})$ (resp. $B(\bm{v})$) is big over $Z$.
\end{enumerate}
Moreover, suppose that one of the following holds:
\begin{itemize}
  \item $K_X+B$ is big over $Z$.
  \item $(X,B)$ is klt and $B$ is big over $Z$.
  \item Conjecture \ref{conj: non-vanishing} holds for $(X/Z,B(\bm{v}))$ for any $\bm{v}:=(v_1,\dots,v_m)\in V_0$.
\end{itemize}
Then for any point $\bm{v}\in V_0$, 
\begin{itemize}
    \item[(4)] $\kappa_{\iota}(X/Z,K_X+B(\bm{v}))=\kappa_{\iota}(X/Z,K_X+B(\bm{v}_0))\geq 0,$ and
    \item[(5)] if either $\bm{v}\in \Qq^c$ or $\kappa(X/Z,K_X+B(\bm{v}_0))\geq 0$, then 
    $$\kappa(X/Z,K_X+B(\bm{v}))=\kappa_{\iota}(X/Z,K_X+B(\bm{v}_0))\geq 0.$$
\end{itemize}
\end{thm}

\begin{proof}
Let $U^1\ni \bm{v}_0$ be an open set satisfies the properties of Corollary \ref{thm: Uniform pe polytope} which only depends on $d$ and $\bm{v}_0$. We will show that the set
$$V_0:=\{\frac{\epsilon'}{\epsilon}\bm{v}_0+(1-\frac{\epsilon'}{\epsilon})\bm{v}\mid \bm{v}\in U^1\}$$
has the required properties. In fact, by our choice, $(X,B(\bm{v}))$ is lc and $K_X+B(\bm{v})$ is pseudo-effective over $Z$ for any $\bm{v}\in U^1$. As
$$K_X+B(\frac{\epsilon'}{\epsilon}\bm{v}_0+(1-\frac{\epsilon'}{\epsilon})\bm{v})=\frac{\epsilon'}{\epsilon}(K_X+B(\bm{v}_0))+(1-\frac{\epsilon'}{\epsilon})(K_X+B(\bm{v})),$$
and $\Supp B(\frac{\epsilon'}{\epsilon}\bm{v}_0+(1-\frac{\epsilon'}{\epsilon})\bm{v})\subset\Supp B(\bm{v}_0)$ for any $\bm{v}\in V_0$, we get (3). Moreover, by the convexity of minimal log discrepancies, 
\begin{align*}
    \mld(X,B(\frac{\epsilon'}{\epsilon}\bm{v}_0+(1-\frac{\epsilon'}{\epsilon})\bm{v}))&\geq\frac{\epsilon'}{\epsilon}\mld(X,B(\bm{v}_0))+(1-\frac{\epsilon'}{\epsilon})\mld(X,B(\bm{v}))
    \geq\frac{\epsilon'}{\epsilon}\mld(X,B(\bm{v}_0)),
\end{align*}
we see that (1) holds. Now we prove (4) and (5). Since $K_X+B(\bm{v})$ is pseudo-effective over $Z$ for any $\bm{v}\in V_0$, by our assumption, $|(K_X+B(\bm{v}))/Z|_{\mathbb R}\not=\emptyset$ for any $\bm{v}\in V_0$. Hence $\kappa_{\iota}(X/Z,K_X+B(\bm{v}))\geq 0$ for any  $\bm{v}\in V_0$. Then (4) follows from Lemma \ref{lem: kod is cons on open set}. 

It suffices to prove (5). If $\bm{v}\in V_0\cap\Qq^c$, then the result follows from (4) and Lemma \ref{lem: kappaiota is kappa for q divisors}. If $\kappa(X/Z,K_X+B(\bm{v}_0))\geq 0$, then by Lemma \ref{lem: rational kod geq 0 and central kod geq 0 imply kod geq 0 on open set},  $\kappa(X/Z,K_X+B(\bm{v}))\geq 0$ for any point $\bm{v}\in V_0$. Then (5) follows from Lemma \ref{lem: kod is cons on open set}.
\end{proof}

\begin{proof}[Proof of Theorem \ref{thm: Uniform kod polytope assume gmm intro}]
The theorem follows from Theorem \ref{thm: Uniform kod polytope assume gmm intro1} immediately.
\end{proof}

\section{Uniform perturbation of DCC coefficients}\label{sec5}

\begin{prop}\label{prop: main prop perturb dcc}
Let $d$ be a positive integer, $\Ii\subset [0,1]$ a DCC set, $\tilde\Ii_0\subset\Ii$ a finite set, and $\epsilon>\epsilon'>0$ two real numbers. Then there exist a finite set $\Ii_0:=\{v^0_1,\dots,v^0_m\}\subset\bar\Ii$, a projection $g: \bar{\Ii}\rightarrow\Ii_0$ (i.e., $g^2=g$), and an open subset $V_0\ni\bm{v}_0:=(v^0_1,\dots,v^0_m)$ of $V$ the rational envelope of $\bm{v}_0$ depending only on $d,\Ii,\tilde\Ii_0,\epsilon$ and $\epsilon'$ satisfying the following. 

Assume that $(X/Z,B:=\sum b_jB^j)$ is an lc (resp. klt, $\epsilon$-lc, $\epsilon$-plt, $\epsilon$-klt) pair of dimension $d$ such that $B\in\Ii$ and each $B^j$ is $\Qq$-Cartier. Then there exist distinct Weil divisors $B_1,\dots, B_m\geq 0$ on $X$ such that
\begin{enumerate}
  \item $g(\gamma)\geq\gamma$ for any $\gamma\in\Ii$ and $g(\gamma)=\gamma$ for any $\gamma\in\tilde\Ii_0$,
  \item $B(\bm{v}_0)=\sum g(b_j)B^j$, where $B(\bm{v}):=\sum_{i=1}^mv_iB_i$ for any $\bm{v}:=(v_1,\dots,v_m)\in\Rr^m$,
  \item $(X,B(\bm{v}))$ is lc for any $\bm{v}\in V_0$,
  \item $(X,B(\bm{v})-D)$ is lc (resp. klt, $\epsilon'$-lc, $\epsilon'$-plt, $\epsilon'$-klt) for any $\bm{v}\in V_0$, where $D:=B(\bm{v}_0)-B$,
  \item if $K_X+B$ is pseudo-effective (resp. big) over $Z$, then $K_X+B(\bm{v})-D$ is pseudo-effective (resp. big) over $Z$ for any $\bm{v}\in V_0$, and
  \item if $B$ is big over $Z$, then $B(\bm{v})-D$ is big over $Z$ for any $\bm{v}\in V_0$.
\end{enumerate}
Moreover, assume that $K_X+B$ is pseudo-effective over $Z$, and one of the following holds.
\begin{itemize}
  \item $K_X+B$ is big over $Z$.
  \item $(X,B)$ is klt and $B$ is big over $Z$.
  \item Conjecture \ref{conj: non-vanishing} holds for $(X/Z,B(\bm{v})-D)$ for any $\bm{v}\in V_0$.
\end{itemize}
Then for any point $\bm{v}\in V_0$,
\begin{itemize}
  \item[(7)] $\kappa_{\iota}(X/Z,K_X+B(\bm{v})-D)=\kappa_{\iota}(X/Z,K_X+B)\geq 0,$ and
  \item[(8)] if either $\bm{v}\in V_0\cap\Qq^m$ or $\kappa(X/Z,K_X+B)\geq 0$, then $$\kappa(X/Z,K_X+B(\bm{v})-D)=\kappa_{\iota}(X/Z,K_X+B)\geq 0.$$
\end{itemize}
\end{prop}

\begin{proof}

\noindent\textbf{Step 1}. In this step, we construct a positive real number $\delta$, finite sets $\Ii_0,\Ii_0',\Ii_0''$, and functions $g,h:\bar\Ii\rightarrow\Ii_0$ where $g$ is a projection.

Take a finite set $\Ii_1\subset [0,1]$ and a function $h_1:\bar\Ii\rightarrow\Ii_1$ depending only on $d$ and $\Ii$ satisfying the properties of Cororllary \ref{cor: gap of psd thresholds}. Since $\Ii$ is a DCC set and $\Ii_1$ is a finite set, one can define
$$\delta:=\min\{\gamma-\gamma_1\mid \gamma\in\bar\Ii,\gamma_1\in\Ii_1,\gamma>\gamma_1\}>0.$$
By \cite[Theorem 5.18]{HLS19}, there exist a finite set $\Ii_2\subset\bar\Ii$ and a projection $g_2: \bar\Ii\rightarrow\Ii_2$ depending only $d$ and $\Ii$, such that $(X,\sum g_2(b_j)B^j)$ is lc, and $\gamma+\frac{1}{3}\delta\geq g_2(\gamma)\geq\gamma$ for any $\gamma\in\bar\Ii$. Let
$$\Ii'_0:=\Ii_2\cup\{\gamma\in\bar\Ii\mid h_1(\gamma)=\gamma\},\ \Ii_0'':=\Ii_1\cup\{\gamma\in\bar\Ii\mid g_2(\gamma)=\gamma\}\text{, and } \Ii_0:=\Ii'_0\cup\Ii''_0.$$
We may define $g$ and $h$ as follows. For any $\gamma\in\bar\Ii$, if $h_1(\gamma)<\gamma<g_2(\gamma)$ and $\gamma\not\in\tilde\Ii_0$, then we let $g(\gamma):=g_2(\gamma)$ and $h(\gamma):=h_1(\gamma)$, otherwise let $g(\gamma):=h(\gamma):=\gamma$.

Define $\bar B:=\sum g(b_j)B^j$ and $\widehat B:=\sum h(b_j)B^j$. By construction, if $K_X+B$ is big over $Z$, then $K_X+\bar B$ and $K_X+\widehat B$ are big over $Z$, and if $B$ is big over $Z$, then $\bar B$ and $\widehat B$ are big over $Z$.

\medskip

\noindent\textbf{Step 2}. In this step, we construct an open set $V_0$ of $V$ and Weil divisors $B_1,\dots,B_m>0$.

Assume that $\Ii_0:=\{v^0_1,\dots,v^0_m\}$ for some real numbers $v^0_1,\dots,v^0_m$. Let $\bm{v}_0:=(v^0_1,\dots,v^0_m)$ and $V\subset\Rr^m$ the rational envelope of $\bm{v}_0$. Then there exist uniquely determined distinct Weil divisors $B_1,\dots,B_m\geq 0$ on $X$, such that $\bar B=\sum_{i=1}^mv^0_1B_i$. We define $B(\bm{v}):=\sum_{i=1}^mv_iB_i$ for any point $\bm{v}:=(v_1,\dots,v_m)\in\Rr^m$, and $D:=B(\bm{v}_0)-B$.

Since $\Ii_0$ is a finite set, there exist a point $\bm{r}:=(r_1,\dots,r_c)$ such that $1,r_1,\dots,r_c$ are linearly independent over $\Qq$, and $\Qq$-linear functions $$s_{1,1},\dots s_{1,p_1},s_{2,1},\dots,s_{2,p_2},s_{3,1},\dots,s_{3,p_3}:\Rr^{c+1}\rightarrow\Rr$$ depending only on $d$ and $\Ii_0$, such that for any $(X,B)$ as in assumption, we may find Weil divisors $$F_{1,1},\dots,F_{1,p_1},F_{2,1},\dots,F_{2,p_2},F_{3,1},\dots,F_{3,p_3}\geq 0$$ on $X$, such that
\begin{itemize}
    \item $\cup_{i=1}^{p_2}\Supp F_{2,i}=\cup_{i=1}^{p_3}\Supp F_{3,i},$
    \item $(\cup_{i=1}^{p_1}\Supp F_{1,i})\cap (\cup_{i=1}^{p_2}\Supp F_{2,i})=\emptyset,$
    \item for any irreducible component $C$ of $\cup_{i=1}^{p_1}\Supp F_{1,i}$, $$h(\mult_CB)=g(\mult_CB)=\mult_CB=s_{1,i}(1,\bm{r}),$$
    \item for any irreducible component $C$ of $\cup_{i=1}^{p_2}\Supp F_{2,i}$, $$s_{3,i}(1,\bm{r})=h(\mult_CB)<\mult_CB<g(\mult_CB)=s_{2,i}(1,\bm{r}),$$
    \item $\bar B=\sum_{i=1}^{p_1}s_{1,i}(1,\bm{r})F_{1,i}+\sum_{i=1}^{p_2}s_{2,i}(1,\bm{r})F_{2,i},$ and
    \item $\widehat B=\sum_{i=1}^{p_1}s_{1,i}(1,\bm{r})F_{1,i}+\sum_{i=1}^{p_3}s_{3,i}(1,\bm{r})F_{3,i}.$
\end{itemize}
For any $\bm{u}\in\Rr^c$, we define 
$$\bar B(\bm{u}):=\sum_{i=1}^{p_1}s_{1,i}(1,\bm{u})F_{1,i}+\sum_{i=1}^{p_2}s_{2,i}(1,\bm{u})F_{2,i}\text{ and }\widehat B(\bm{u}):=\sum_{i=1}^{p_1}s_{1,i}(1,\bm{u})F_{1,i}+\sum_{i=1}^{p_3}s_{3,i}(1,\bm{u})F_{3,i}.$$
Since $s_{i,j}(1\le i\le 3,1\le j\le p_i)$ are continuous functions,
\begin{itemize}
    \item by \cite[Corollary 5.5]{HLS19}, there exists a convex open subset $U_1\ni\bm{r}$ of $\Rr^c$ depending only on $d$ and $s_{1,1},\dots,s_{1,p_1}$ and $s_{2,1},\dots,s_{2,p_2}$, such that $(X,\bar B(\bm{u}))$ is lc, and $||\bar B-\bar B(\bm{u})||<\frac{\delta}{3}$, for any $\bm{u}\in U_1$, and
    \item by Theorem \ref{thm: one variable psd polytope}, there exists a convex open subset $U_2\ni\bm{r}$ of $\Rr^c$ depending only on $d$, $s_{1,i},\dots,s_{1,p_1}$ and $s_{3,1},\dots,s_{3,p_3}$, such that $(X,\widehat B(\bm{u}))$ is lc and $K_X+\widehat B(\bm{u})$ is pseudo-effective over $Z$ for any $\bm{u}\in U_2$.
\end{itemize}
We define 
$$U_0:=\{\frac{\epsilon'}{\epsilon}\bm{r}+(1-\frac{\epsilon'}{\epsilon})\bm{u}\mid \bm{u}\in U_1\cap U_2\}.$$

By construction, there exist uniquely determined distinct Weil divisors $B_1',\dots, B_m'\geq 0$ on $X$, such that $\widehat B=\sum_{i=1}^mv_iB_i'$. For any point $\bm{v}:=(v_1,\dots,v_m)\in\Rr^m$, we define $B'(\bm{v}):=\sum_{i=1}^mv_iB_i'$. Moreover, there exists a $\Qq$-affine function $l: \Rr^c\rightarrow V$, such that
\begin{itemize}
    \item for any $\bm{u}\in\Rr^c$, $\bar B(\bm{u})=B(l(\bm{u}))$ and $\widehat B(\bm{u})=B'(l(\bm{u}))$, and
    \item $l(\bm{r})=\bm{v}_0$, in particular, $\bar B(\bm{r})=B(\bm{v}_0)$ and $\widehat B(\bm{r})=B'(\bm{v}_0)$.
\end{itemize}
We let $V_0:=l(U_0)\cap V(d,\epsilon,\epsilon',\bm{v}_0)$, where $V(d,\epsilon,\epsilon',\bm{v}_0)$ is an open set satisfying the properties of Theorem \ref{thm: Uniform kod polytope assume gmm intro1} which only depends on $d,\epsilon,\epsilon',$ and $\bm{v}_0$.

\medskip

\noindent\textbf{Step 3}. In this step, we prove (1--6). 

It is obvious that (1) and (2) follow from our construction, and (3) follows from Claim \ref{claim: xbbmv is lc for every bmv in u'}.

\begin{claim}\label{claim: xbbmv is lc for every bmv in u'}
$(X,B(\bm{v}))$ is lc for every $\bm{v}\in V'$, where $V':=l(U_1\cap U_2)\supset V_0$.
\end{claim}
\begin{proof}[Proof of Claim \ref{claim: xbbmv is lc for every bmv in u'}]
For every $\bm{v}\in V'$, there exists $\bm{u}\in U_1\cap U_2$ such that $B(\bm{v})=\bar B(\bm{u}).$ Since $\bm{u}\in U_1\cap U_2\subset U_1$, $(X,\bar B(\bm{u}))$ is lc and hence $(X,B(\bm{v}))$ is lc.
\end{proof}

\begin{claim}\label{claim: in u1capu2 bbmu-dgeqtilde bbmugeq 0}
Fix any $\bm{u}\in U_1\cap U_2$. For any prime divisor $C$,
\begin{enumerate}
    \item[(i)] if $C\subset \Supp D$, then $\mult_C(\bar B(\bm{u})-D)>\mult_C\widehat B(\bm{u})\geq 0,$ and
    \item[(ii)] if $C\subset\Supp(\widehat B(\bm{u})+\bar B(\bm{u}))$ and $C\nsubseteq \Supp D$, then $\mult_C\bar B(\bm{u})=\mult_C\widehat B(\bm{u})\ge0$.
\end{enumerate}
In particular, $\bar B(\bm{u})-D\geq\widehat B(\bm{u})\geq 0$.
\end{claim}
\begin{proof}[Proof of Claim \ref{claim: in u1capu2 bbmu-dgeqtilde bbmugeq 0}]
Recall that $D:=B(\bm{v}_0)-B=\sum (g(b_j)-b_j)B^j$. Thus
$$\Supp D=\cup_{\{j\mid g(b_j)>b_j\}}\Supp B^j=\cup_{i=1}^{p_2}\Supp F_{2,i}=\cup_{i=1}^{p_3}\Supp F_{3,i}.$$
By the construction of $\bar B(\bm{u})$ and $\widehat B(\bm{u})$, we have (ii). Since $(X,\widehat B(\bm{u}))$ is lc, $\widehat B(\bm{u})\geq 0$. For any irreducible component $C$ of $D$, the following hold.
\begin{itemize}
  \item Since $0\leq g(b_j)-b_j\leq\frac{1}{3}\delta$ for any $j$, $0<\mult_CD\leq\frac{1}{3}\delta$.
  \item As $C\subset\cup_{i=1}^{p_2}\Supp F_{2,i}$, $h(\mult_CB)<\mult_CB<g(\mult_CB).$
  \item By the construction of $\delta$, $\mult_C(B-\widehat B)\geq\delta$.
\end{itemize}
Thus 
\begin{align*}
    \mult_C(\bar B(\bm{u})-D)=&\mult_C(\bar B(\bm{u})-\bar B)+\mult_C(\bar B-B)+\mult_C(B-\widehat B)\\
    &+\mult_C(\widehat B-\widehat B(\bm{u}))+\mult_C(\widehat B(\bm{u}))-\mult_CD\\
    >&-\frac{\delta}{3}+0+\delta-\frac{\delta}{3}+\mult_C(\widehat B(\bm{u}))-\frac{\delta}{3}
    =\mult_C(\widehat B(\bm{u}))\geq 0,
\end{align*}
which implies (i).
\end{proof}

\noindent\textit{Proof of Proposition \ref{prop: main prop perturb dcc} continued}. 
For any point $\bm{v}\in V_0$, we may pick $\bm{v}'\in V'$ such that
$$\bm{v}=\frac{\epsilon'}{\epsilon}\bm{v}_0+(1-\frac{\epsilon'}{\epsilon})\bm{v}'.$$
By Claim \ref{claim: xbbmv is lc for every bmv in u'}, $(X,B(\bm{v}'))$ is lc, thus $(X,B(\bm{v}')-D)$ is sub-lc. Then by Claim \ref{claim: in u1capu2 bbmu-dgeqtilde bbmugeq 0}, $(X,B(\bm{v}')-D)$ is lc. Since $(X,B)$ is lc (resp. klt, $\epsilon$-lc, $\epsilon$-plt, $\epsilon$-klt) and $B=B(\bm{v}_0)-D$, $(X,B(\bm{v})-D)$ is lc (resp. klt, $\epsilon'$-lc, $\epsilon'$-plt, $\epsilon'$-klt), which implies (4).

Under the assumption of (5) (resp. (6)), $K_X+\widehat B$ is pseudo-effective (resp. big, $\widehat B$ is big) over $Z$. Thus $K_X+\widehat B(\bm{u})$ is pseudo-effective (resp. big, $B(\bm{u})$ is big) over $Z$ for any point $\bm{u}\in U_0$. By Claim \ref{claim: in u1capu2 bbmu-dgeqtilde bbmugeq 0},  $\bar B(\bm{u})-D-\widehat B(\bm{u})\geq 0$ for any point $\bm{u}\in U_0$, which implies (5) and (6).

\medskip

\noindent\textbf{Step 4}. In this step, we prove the ``Moreover'' part, and hence finish the proof. 

Indeed, by (5), (6) and \cite[Theorem D]{BCHM10}, $\kappa_{\iota}(X/Z,K_X+B(\bm{v})-D)\geq 0$ for every $\bm{v}\in V_0$. Then (7) follows from Lemma \ref{lem: kod is cons on open set}. It suffices to prove (8). 
\begin{claim}\label{claim: barbu-d rational decomposition}
For any point $\bm{u}\in U_1\cap U_2$, there exists a rational number $t\in (0,1)$ such that
$$\bar B(\bm{u})-D\geq t\bar B(\bm{u})+(1-t)\widehat B(\bm{u}).$$
\end{claim}
\begin{proof}[Proof of Claim \ref{claim: barbu-d rational decomposition}]
If $D=0$, then we take $t:=\frac{1}{2}$. From now on we assume that $D\not=0$. Let
$$\beta:=\min\{\mult_C(\bar B(\bm{u})-D-\widehat B(\bm{u}))\mid C\text{ is an irreducible component of }D\}.$$
By Claim \ref{claim: in u1capu2 bbmu-dgeqtilde bbmugeq 0}(i), $\beta>0$. Since $(X,\bar B(\bm{u}))$ is lc, $\beta<1$. Choose any rational number $t\in (0,\beta)$. Pick a prime divisor $C$ on $X$. If $C\subseteq\Supp D$, then
\begin{align*}
    &\ \ \ \mult_C(\bar B(\bm{u})-D)=\mult_C(\bar B(\bm{u})-D-\widehat B(\bm{u}))+\mult_C\widehat B(\bm{u})\geq t+\mult_C\widehat B(\bm{u})
    \\&\geq t(\mult_C(\bar B(\bm{u})-\widehat B(\bm{u})))+\mult_C\widehat B(\bm{u})
    =\mult_C(t\bar B(\bm{u})+(1-t)\widehat B(\bm{u})).
\end{align*}
If $C\subset\Supp(\bar B(\bm{u})+\widehat B(\bm{u}))$ and $C\nsubseteq\Supp D$, then by Claim \ref{claim: in u1capu2 bbmu-dgeqtilde bbmugeq 0}(ii), we have
$$\mult_C(\bar B(\bm{u})-D)=\mult_C(t\bar B(\bm{u})+(1-t)\widehat B(\bm{u})).$$
Thus $t$ has the required property.
\end{proof}

\begin{claim}\label{claim: barbu-d rational decomposition perturbed}
For any point $\bm{u}\in l^{-1}(V_0)$, there exists a rational number $t\in (0,1)$, a point $\bm{y}\in l^{-1}(V_0)$, and a point $\bm{z}\in l^{-1}(V_0)\cap\Qq^c$, such that
$$\bar B(\bm{u})-D\geq t\bar B(\bm{y})+(1-t)\widehat B(\bm{z}).$$
\end{claim}
\begin{proof}[Proof of Claim \ref{claim: barbu-d rational decomposition perturbed}]
If $D=0$, then we pick $t:=\frac{1}{2}$. Since $l^{-1}(V_0)$ is an open set, we may find $\bm{z}\in l^{-1}(V_0)\cap\Qq^c$ such that $\bm{y}:=\bm{u}+(\bm{u}-\bm{z})\in l^{-1}(U)$. Since $\widehat B(\bm{z})=\bar B(\bm{z})$ in this case, the claim holds. From now on we assume that $D\not=0$.  Let
$$\beta:=\min\{\mult_C(\bar B(\bm{u})-D-\widehat B(\bm{u}))\mid C\text{ is an irreducible component of }D\}.$$
By Claim \ref{claim: in u1capu2 bbmu-dgeqtilde bbmugeq 0}(i), $\beta>0$. Since $(X,\bar B(\bm{u}))$ is lc, $\beta<1$. Let $t\in (0,\beta)$ be any rational number. Then there exists a positive real number $\alpha$, such that for any $\bm{x}\in\Rr^c$, if $||\bm{x}-\bm{u}||_{\infty}<\alpha$, then $||\bar B(\bm{x})-\bar B(\bm{u})||<\beta-t$, $||\widehat B(\bm{x})-\widehat B(\bm{u})||<\beta-t$, and $\bm{x}\in l^{-1}(V_0)$.

Let $C$ be a prime divisor on $X$. Pick any $\bm{z}\in\Qq^m$ such that $||\bm{z}-\bm{u}||_{\infty}<t\alpha$, and let $$\bm{y}:=\bm{z}+\frac{1}{t}(\bm{u}-\bm{z}).$$ 
Then $||\bm{y}-\bm{u}||<(1-t)\alpha$. Thus if $C\subset\Supp D$, then
\begin{align*}
    &\mult_C(\bar B(\bm{u})-D)=\mult_C(\bar B(\bm{u})-D-\widehat B(\bm{u}))+\mult_C\widehat B(\bm{u})
    \geq \beta+\mult_C\widehat B(\bm{u})\\
    =& (\beta-t)+t+\mult_C\widehat B(\bm{u})
    \geq(\beta-t)+t(\mult_C(\bar B(\bm{u})-\widehat B(\bm{u})))+\mult_C\widehat B(\bm{u})\\
    =&(\beta-t)+\mult_C(t\bar B(\bm{u})+(1-t)\widehat B(\bm{u}))\\
    =& (\beta-t)+\mult_C(t\bar B(\bm{y})+(1-t)\widehat B(\bm{z}))
    +t\mult_C(\bar B(\bm{u})-\bar B(\bm{y}))+(1-t)\mult_c(\bar B(\bm{u})-\bar B(\bm{z}))\\
    >&(\beta-t)+\mult_C(t\bar B(\bm{y})+(1-t)\widehat B(\bm{z}))
    -t(\beta-t)-(1-t)(\beta-t)\\
    =&\mult_C(t\bar B(\bm{y})+(1-t)\widehat B(\bm{z})).
\end{align*}
If $C\subset\Supp(\bar B(\bm{u})+\widehat B(\bm{u}))$ and $C\nsubseteq\Supp D$, then by the construction of $\bar B(\bm{u}),\bar B(\bm{y})$, $\widehat B(\bm{z})$, and the fact that $\bm{u}=t\bm{y}+(1-t)\bm{z}$, we see that
$$\mult_{C}(\bar B(\bm{u})-D)=\mult_C\bar B(\bm{u})=\mult_C(t\bar B(\bm{y})+(1-t)\widehat B(\bm{z})).$$
The claim follows.
\end{proof}

\noindent\textit{Proof of Proposition \ref{prop: main prop perturb dcc} continued}. 
By the assumption, $\kappa_{\iota}(X/Z,K_X+\widehat B)\geq 0$ which implies that $\kappa_{\iota}(X/Z,K_X+\bar B)\geq 0$. By the choices of $V(d,\epsilon,\epsilon',\bm{v}_0)\supset V_0$, $\kappa_{\iota}(X/Z,K_X+\widehat B(\bm{u}))\geq 0$ for any $\bm{u}\in l^{-1}(V_0)$. By Claim \ref{claim: in u1capu2 bbmu-dgeqtilde bbmugeq 0}, $\bar B(\bm{u})\geq\widehat B(\bm{u})$ for any $\bm{u}\in l^{-1}(V_0)$, so $\kappa_{\iota}(X/Z,K_X+\bar B(\bm{u}))\geq 0$. By Claim \ref{claim: barbu-d rational decomposition}, for any point $\bm{u}\in l^{-1}(V_0)$, there exists a rational number $t\in (0,1)$ such that 
\begin{equation}\label{eq1}
B(\bm{v})-D=\bar B(\bm{u})-D\geq t\bar B(\bm{u})+(1-t)\widehat B(\bm{u}).
\end{equation}

Suppose that $\bm{v}\in V_0\cap\Qq^m$. We may pick $\bm{u}\in l^{-1}(\bm{v})\cap\Qq^c\subset l^{-1}(V_0)\cap\Qq^c$. Then by Lemma \ref{lem: kappaiota is kappa for q divisors}, $\kappa(X/Z,K_X+\bar B(\bm{u}))\geq 0$ and $\kappa(X/Z,K_X+\widehat B(\bm{u}))\geq 0$, which imply that $\kappa(X/Z,K_X+B(\bm{v})-D)\geq 0$ by \eqref{eq1}. By (7) and Proposition \ref{prop: basic properties of three iitaka dimensions}, we have
$$\kappa(X/Z,K_X+B(\bm{v})-D)=\kappa_{\iota}(X/Z,K_X+B)\geq 0.$$

Now assume that $\kappa(X/Z,K_X+B)\geq 0$. Since $\bar B\geq B$, $\kappa(X/Z,K_X+\bar B)\geq 0$. For any $\bm{v}\in V_0$, pick $\bm{u}\in l^{-1}(\bm{v})\subset l^{-1}(V_0)$. By Claim \ref{claim: barbu-d rational decomposition perturbed}, there exist a rational number $t\in (0,1)$, $\bm{y}\in l^{-1}(V_0)$, and $\bm{z}\in l^{-1}(V_0)\cap\Qq^c$, such that
\begin{equation}\label{eq2}
    B(\bm{v})-D=\bar B(\bm{u})-D\geq t\bar B(\bm{y})+(1-t)\widehat B(\bm{z}).
\end{equation}
Since $l(\bm{y})\in V\subset V(d,\epsilon,\epsilon',\bm{v}_0)$, by Theorem \ref{thm: Uniform kod polytope assume gmm intro1}(5), \begin{equation}\label{eq3}
\kappa(X/Z,K_X+\bar B(\bm{y}))=\kappa(X/Z,K_X+B(l(\bm{y})))=\kappa_{\iota}(X/Z,K_X+\bar{B})\geq 0.
\end{equation}
Since $\bm{z}\in l^{-1}(V_0)\cap\Qq^c$, by Proposition \ref{prop: basic properties of three iitaka dimensions}, one can see that
\begin{equation}\label{eq4}
\kappa(X/Z,K_X+\widehat B(\bm{z}))\geq 0.
\end{equation}
Thus $\kappa(X/Z,K_X+B(\bm{v})-D)\geq 0$ by \eqref{eq2}, \eqref{eq3} and \eqref{eq4}. Then by (7) and Proposition \ref{prop: basic properties of three iitaka dimensions}, 
$$\kappa(X/Z,K_X+B(\bm{v})-D)=\kappa_{\iota}(X/Z,K_X+B)\geq 0,$$
and thus (8) holds.
\end{proof}

\begin{thm}\label{thm:dccdecomp1}
Let $d$ be a positive integer, $\Ii\subset[0,1]$ a DCC set, and $\epsilon>\epsilon'>0$ two real numbers. Then there exist a finite set $\Ii'\subset[0,1]\cap\Qq$, an ACC set $\Ii''\subset[0,1]$, and real numbers $a_1,\dots,a_k\in (0,1]$ depending only on $d$, $\Ii$, $\epsilon$ and $\epsilon'$ satisfying the following. 

Assume that $(X/Z,B)$ is a $d$-dimensional lc (resp. klt, $\epsilon$-lc, $\epsilon$-plt, $\epsilon$-klt) pair such that $B\in\Ii$, all the irreducible components of $B$ are $\Qq$-Cartier, and $K_X+B$ is pseudo-effective over $Z$. Then there exist an $\Rr$-divisor $D\geq 0$ and $\Qq$-divisors $B_1,\dots,B_k\geq 0$ on $X$, such that 
\begin{enumerate}
    \item $\sum_{i=1}^ka_i=1$,
    \item $\sum_{i=1}^ka_iB_i-D=B$,
    \item $B_i\in\Ii'$ for each $i$, and $D\in\Ii''$,
    \item $(X,B_i)$ and $(X,B_i-D)$ are lc (resp. klt, $\epsilon'$-lc, $\epsilon'$-plt,
    $\epsilon'$-klt) for each $i$, 
    \item $K_X+B_i-D$ is pseudo-effective over $Z$ for each $i$, and
    \item if $K_X+B$ (resp. $B$) is big over $Z$, then $K_X+B_i-D$ (resp. $B_i-D$) is big ovr $Z$ for each $i$.
\end{enumerate}
Moreover, suppose that one of the following holds:
\begin{itemize}
  \item $K_X+B$ is big over $Z$.
  \item $(X,B)$ is klt and $B$ is big over $Z$.
  \item Conjecture \ref{conj: non-vanishing} holds.
\end{itemize}
Then $\kappa(X/Z,K_X+B_i-D)=\kappa_{\iota}(X/Z,K_X+B)\geq 0$ for each $i$.
\end{thm}

\begin{proof}
We may pick a point $\bm{v}_0:=(v^0_1,\dots,v^0_m)\in\Rr^m$ and an open subset $V_0$ which only depends on $d,\Ii,\epsilon$, and $\epsilon'$ satisfying the properties of Proposition \ref{prop: main prop perturb dcc}. For any $(X,B)$ as in the assumption, let $B(\bm{v})$ be the $\Rr$-divisor as in Proposition \ref{prop: main prop perturb dcc} and $D:=B(\bm{v}_0)-B$. It is clear that $D\in\Ii''$ for some ACC set $\Ii''\subset [0,1]$ which only depends on $d,\Ii,\epsilon$, and $\epsilon'$.

Let $k:=\dim V_0+1$ and $\bm{v}_1,\dots,\bm{v}_k\in V_0\cap\Qq^m$ points depending only on $d,\Ii,\epsilon$, and $\epsilon'$, such that $\bm{v}_0$ is contained in the interior of the convex hull spanned by $\bm{v}_1,\dots,\bm{v}_k$. Thus we may find real numbers $a_1,\dots,a_k\in (0,1]$, such that $\sum_{i=1}^ka_i=1$ and $\sum_{i=1}^ka_i\bm{v}_i=\bm{v}_0$. Set $B_i:=B(\bm{v}_i)$ for each $i$. Then $B_i\in\Ii$ for some finite set $\Ii\subset\Qq\cap[0,1]$ depending only on $d,\Ii,\epsilon$, and $\epsilon'$. Now (1--3) follow from our construction, and (4--8) follow from Proposition \ref{prop: main prop perturb dcc}.
\end{proof}
\begin{proof}[Proof of Theorem \ref{thm:dccdecomp}]
According to Theorem \ref{thm:dccdecomp1}, there exist a finite set $\Ii'\subset[0,1]\cap\Qq$, an ACC set $\Ii''\subset[0,1]$, and real numbers $a_1,\dots,a_k\in (0,1]$ depending only on $d$ and $\Ii$ satisfying the properties of Theorem \ref{thm:dccdecomp1}. Let
$$\Ii_0:=\{a-b\mid a\in\Ii',b\in\Ii''\}\cap[0,1].$$
It is clear that $a_1,\dots,a_k$ and $\Ii_0$ have the required properties.
\end{proof}

\begin{rem}
Theorems \ref{thm: Uniform kod polytope assume gmm intro1} and \ref{thm:dccdecomp1} can be generalized to generalized pairs (g-pairs) easily by using very similar arguments and the tools developed in \cite{BZ16,HL21a,HLi18}. In fact, the structure of generalized pairs will naturally appear when we study (effective) Iitaka fibrations, as whenever we have an lc-trivial fibration structure, generalized pairs with DCC coefficients can be found in the canonical bundle formulas. Although we have a decomposable canonical bundle formula for generalized pairs with finite real coefficients (cf. \cite{HL21b,HJL22,JLX22}), decomposable canonical bundle formulas for generalized pairs with arbitrary DCC coefficients remain to be studied.  We will discuss this in details in our sequel paper.
\end{rem}

\end{document}